\documentclass[reqno]{amsart}
\usepackage{hyperref}

\begin{document}
\title[\hfilneg Non-existence of global solutions]
{Non-existence of global solutions for a generalized fractional differential problem}

\author[Sandeep P Bhairat]
{Sandeep P Bhairat}

\address{Sandeep P Bhairat \newline
Department of Mathematics,
Institute of Chemical Technology,\newline
Mumbai -- 400 019, (M.S) India.}
\email{sp.bhairat@ictmumbai.edu.in}


\keywords{Non-existence; Global solution; Fractional differential equations;
\hfill\break\indent Katugampola fractional integral and derivative; Test function method}

\begin{abstract}
Aim of this paper is to study the non-existence of global solutions of the fractional differential problem involving generalized Katugampola derivative. We utilize the test function method and fractional integration by parts formula to obtain the result. An illustrative example is also given.
\end{abstract}

\maketitle
\numberwithin{equation}{section}
\newtheorem{theorem}{Theorem}[section]
\newtheorem{lemma}[theorem]{Lemma}
\newtheorem{remark}[theorem]{Remark}
\newtheorem{definition}[theorem]{Definition}
\newtheorem{example}[theorem]{Example}
\allowdisplaybreaks

\section{Introduction}
In the past thirty years, the interest to fractional differential equations paid more attention of many researchers in several areas such as bioengineering, physics, mechanics and applied sciences, \cite{hr2,ve}. For recent development and historical arguments, see the monographs \cite{kst,skm}. The existence of solutions for various class of fractional differential equations are studied extensively with number of fractional derivatives in \cite{as}-\cite{kmf},\cite{mdk,ke,oo}. Whereas for the non-existence of solutions, one can see the recent papers \cite{neh,nehh,newnec,newne,oldnex}.

Recently, the existence and uniqueness of solution of initial value problem 
\begin{equation}\label{p}
\begin{cases}
& ^{\rho}D_{a+}^{\alpha,\beta}x(t)=f(t,x(t)),\quad 0<\alpha<1,0\leq\beta\leq1,
\rho>0, \, t>a>0,\\
& ^{\rho}I_{a+}^{1-\gamma }x(t) \big|_{t=a}=\phi,\quad \qquad\phi\in \mathbb{R},\gamma=\alpha-\beta-\alpha\beta,
\end{cases}
\end{equation}
for the generalized fractional differential problem is studied in \cite{oo}.  The operators $^{\rho}I_{a+}^{\sigma }$ and $^{\rho}D_{a+}^{\sigma,\xi }$ are the Katugampola fractional integral \cite{ki} and generalized Katugampola fractional derivative \cite{oo}, respectively. In this paper we consider the generalized fractional differential problem of type
\begin{equation}\label{pex}
\begin{cases}
& ^{\rho}D_{a+}^{\alpha,\beta}x(t) \geq {\big(\frac{t^\rho-a^\rho}{\rho}\big)}^{\mu}{|x(t)|}^{m},\quad \rho>0,\, t>a>0,\,\mu\in\mathbb{R},\,\,m>1,\\
& ^{\rho}I_{a+}^{1-\gamma }x(t) \big|_{t=a}=x_a,\qquad\qquad x_a\in \mathbb{R},\,\gamma=\alpha+\beta-\alpha\beta,
\end{cases}
\end{equation}
where $0<\alpha<1,0\leq \beta\leq1$. We prove that no solutions can exist for all time for certain values of $\mu$ and $m$ in an appropriate weighted space of continuous functions.
\section{Preliminaries}
In this section, we list some definitions and lemmas useful throughout the paper.
\begin{definition}\label{d1}\cite{kst}
The space $X_{c}^{p}(a,b)\,(c\in\mathbb{R},p\geq1)$ consists of those real-valued Lebesgue measurable functions $g$ on $(a,b)$ for which ${\|g\|}_{X_{c}^{p}}<\infty,$ where
\begin{gather*}
{\|g\|}_{X_{c}^{p}}={\bigg(\int_{a}^{b}{|t^cg(t)|}^{p}\frac{dt}{t}\bigg)}^{\frac{1}{p}},\quad p\geq1,\,\,c\in\mathbb{R},\\
{\|g\|}_{X_{c}^{p=\infty}}=\text{ess sup}_{a\leq t\leq b}|t^cg(t)|,\quad c\in\mathbb{R}.
\end{gather*}
In particular, when $c=\frac{1}{p},$ we see that $X_{{1}/{p}}^{c}(a,b)=L_p(a,b).$
\end{definition}
\begin{definition}\label{d2}\cite{oo}
Let $\Omega=[a,b]\,\,(0<a<b<\infty)$ be a finite interval on $\mathbb{R}^{+}$ and $\rho>0.$ Denote by $C[a,b]$ a space of continuous functions $g$ on $\Omega$ with the norm
\begin{gather*}
{\|g\|}_{C}=\max_{t\in\Omega}|g(t)|.
\end{gather*}
The weighted space $C_{\gamma,\rho}[a,b]$ of functions $g$ on $(a,b]$ is defined by
\begin{gather}\label{space}
C_{\gamma,\rho}[a,b]=\{g:(a,b]\to\mathbb{R}:{\bigg(\frac{t^\rho-a^\rho}{\rho}\bigg)}^{\gamma}g(t)\in{C[a,b]}\},\quad0\leq\gamma<1
\end{gather}
with the norm
\begin{gather*}
{\|g\|}_{C_{\gamma,\rho}}={\bigg\|{\bigg(\frac{t^\rho-a^\rho}{\rho}\bigg)}^{\gamma}g(t)\bigg\|}_{C}=\max_{t\in\Omega}\bigg|{\bigg(\frac{t^\rho-a^\rho}{\rho}\bigg)}^{\gamma}g(t)\bigg|,
\end{gather*}
and $C_{0,\rho}[a,b]=C[a,b].$
\end{definition}
\begin{definition}\label{d3}\cite{oo}
Let $\delta_\rho=\big(t^{\rho-1}\frac{d}{dt}\big),\,\Omega=[a,b]\,(0<a<b<\infty),\rho>0$ and $0\leq\gamma<1.$ Denote $C_{\delta_\rho,\gamma}^{n}[a,b]$ the Banach space of functions $g$ which are continuously differentiable, with $\delta_\rho,$ on $[a,b]$ upto $(n-1)$ order and have the derivative $\delta_\rho^ng$ on $(a,b]$ such that $\delta_\rho^ng\in{C_{\gamma,\rho}[a,b]}:$
\begin{gather*}
C_{\delta_\rho,\gamma}^{n}[a,b]=\big\{\delta_\rho^kg\in{C[a,b]}, k=0,1,\cdots,n-1,\,\,\delta_\rho^ng\in{C_{\gamma,\rho}[a,b]}\big\},\quad n\in\mathbb{N},
\end{gather*}
with the norm
\begin{gather*}
{\|g\|}_{C_{\delta_\rho,\gamma}^{n}}=\sum_{k=0}^{n-1}{\|\delta_\rho^kg\|}_{C}+{\|\delta_\rho^ng\|}_{C_{\gamma,\rho}},\quad
{\|g\|}_{C_{\delta_\rho}^{n}}=\sum_{k=0}^{n}\max_{t\in\Omega}|\delta_\rho^kg(t)|.
\end{gather*}
Note that, for $n=0,$ we have $C_{\delta_\rho,\gamma}^{0}[a,b]=C_{\gamma,\rho}[a,b].$
\end{definition}
\begin{definition}\label{d4}\cite{ki}
Let $\alpha>0$ and $h\in{X_{c}^{p}(a,b)},$ where $X_{c}^{p}$ is as in Definition \ref{d1}. The left-sided Katugampola fractional integral $^{\rho}I_{a+}^{\alpha}$  of order $\alpha$ is defined by
\begin{gather}\label{kil}
^{\rho}I_{a+}^{\alpha}h(t)=\int_{a}^{t}s^{\rho-1}{\bigg(\frac{t^\rho-s^\rho}{\rho}\bigg)}^{\alpha-1}\frac{h(s)}{\Gamma(\alpha)}ds,\quad t>a.
\end{gather}
\end{definition}
\begin{definition}\label{d5}\cite{ki}
Let $\alpha>0$ and $h\in{X_{c}^{p}(a,b)},$ where $X_{c}^{p}$ is as in Definition \ref{d1}. The right-sided Katugampola fractional integral $^{\rho}I_{b-}^{\alpha}$  of order $\alpha$ is defined by
\begin{gather}\label{kir}
^{\rho}I_{b-}^{\alpha}h(t)=\int_{t}^{b}s^{\rho-1}{\bigg(\frac{s^\rho-t^\rho}{\rho}\bigg)}^{\alpha-1}\frac{h(s)}{\Gamma(\alpha)}ds,\quad t<b.
\end{gather}
\end{definition}
\begin{definition}\label{d6}\cite{kd}
Let $\alpha\in{\mathbb{R}^{+}{\setminus}\mathbb{N}}$ and $n=[\alpha]+1,$ where $[\alpha]$ is integer part of $\alpha$ and $\rho>0.$ The left-sided Katugampola fractional derivative $^{\rho}D_{a+}^{\alpha}$ is defined by
\begin{align}\label{kdl}
^{\rho}D_{a+}^{\alpha}h(t)={\bigg(t^{\rho-1}\frac{d}{dt}\bigg)}^{n}\int_{a}^{t}s^{\rho-1}{\bigg(\frac{t^\rho-s^\rho}{\rho}\bigg)}^{n-\alpha-1}\frac{h(s)}{\Gamma(n-\alpha)}ds.
\end{align}
\end{definition}
\begin{definition}\label{d7}\cite{kd}
Let $\alpha\in{\mathbb{R}^{+}\setminus\mathbb{N}}$ and $n=[\alpha]+1,$ where $[\alpha]$ is integer part of $\alpha$ and $\rho>0.$ The right-sided Katugampola fractional derivative $^{\rho}D_{b-}^{\alpha}$ is defined by
\begin{align}\label{kdr}
^{\rho}D_{b-}^{\alpha}h(t)={\bigg(-t^{\rho-1}\frac{d}{dt}\bigg)}^{n}\int_{t}^{b}s^{\rho-1}{\bigg(\frac{s^\rho-t^\rho}{\rho}\bigg)}^{n-\alpha-1}\frac{h(s)}{\Gamma(n-\alpha)}ds.
\end{align}
\end{definition}
\begin{lemma}\label{l1}\cite{ki}
Suppose that $\alpha>0,\beta>0,p\geq1,0<a<b<\infty$ and $\rho,c\in\mathbb{R}$ such that $\rho\geq{c}.$ Then, for $h\in{X_{c}^{p}(a,b)},$ the semigroup property of Katugampola integral is valid. This is
\begin{gather}\label{ski}
^{\rho}I_{a+}^{\alpha}{^{\rho}I_{a+}^{\beta}h(t)}={^{\rho}I_{a+}^{\alpha+\beta}h(t)}.
\end{gather}
\end{lemma}
A similar property for right-sided operator also holds.
\begin{lemma}\label{l2}\cite{kd}
Suppose that ${^{\rho}I_{a+}^{\alpha}}$ and $^{\rho}D_{a+}^{\alpha}$ are as defined in Definitions \ref{d4} and \ref{d6}, respectively. Then,
\begin{description}
  \item[(i)] ${^{\rho}I_{a+}^{\alpha}}{\big(\frac{t^\rho-a^\rho}{\rho}\big)}^{\sigma-1}=\frac{\Gamma(\sigma)}{\Gamma(\sigma+\alpha)}{\big(\frac{t^\rho-a^\rho}{\rho}\big)}^{\alpha+\sigma-1},\qquad\alpha\geq0,\sigma>0, t>a.$
  \item[(ii)] ${^{\rho}D_{a+}^{\alpha}}{\big(\frac{t^\rho-a^\rho}{\rho}\big)}^{\alpha-1}=0,\qquad0<\alpha<1, t>a.$
\end{description}
\end{lemma}
\begin{lemma}\label{l3}\cite{oo}
Let $0<a<b<\infty,\alpha>0$ and $0\leq\gamma<1.$\\
(i) If $0<\alpha<\gamma,$ then $^{\rho}I_{a+}^{\alpha}$ is bounded from $C_{\gamma,\rho}[a,b]$ into $C_{\gamma,\rho}[a,b].$\\
(ii) If $\gamma\leq\alpha$, then $^{\rho}I_{a+}^{\alpha}$ is bounded from $C_{\gamma,\rho}[a,b]$ into $C[a,b].$\\
In particular, $^{\rho}I_{a+}^{\alpha}$ is bounded in $C_{\gamma,\rho}[a,b].$
\end{lemma}
\begin{lemma}\label{l4}
Let $0\leq\alpha<1$ and $0\leq\gamma<1.$ If $h\in{C_{\gamma,\rho}^1[a,b]},$ then the fractional derivatives $^{\rho}D_{a+}^{\alpha}$ and $^{\rho}D_{b-}^{\alpha}$ exist on $(a,b]$ and $[a,b),$ respectively $(a>0)$ and are represented in the forms:
\begin{gather}\label{kd1}
  ^{\rho}D_{a+}^{\alpha}h(t)=\frac{h(a)}{\Gamma(1-\alpha)}{\bigg(\frac{t^\rho-a^\rho}{\rho}\bigg)}^{-\alpha}+\int_{a}^{t}s^{\rho-1}{\bigg(\frac{t^\rho-s^\rho}{\rho}\bigg)}^{-\alpha}\frac{h^{'}(s)}{\Gamma(1-\alpha)}ds,\\
^{\rho}D_{b-}^{\alpha}h(t)=\frac{h(b)}{\Gamma(1-\alpha)}{\bigg(\frac{b^\rho-t^\rho}{\rho}\bigg)}^{-\alpha}-\int_{t}^{b}s^{\rho-1}{\bigg(\frac{s^\rho-t^\rho}{\rho}\bigg)}^{-\alpha}\frac{h^{'}(s)}{\Gamma(1-\alpha)}ds.
\end{gather}\label{kd2}
\end{lemma}
\begin{lemma}\label{l5}[Fractional integration by parts]
If $\alpha>0,p\geq1,0\leq a<b\leq\infty$ and $\rho,c\in\mathbb{R}$ be such that $\rho\geq{c},$ then, for $g,h\in{X_{c}^{p}(a,b)},$ the following relation hold:
\begin{gather}\label{fip}
\int_{a}^{b}t^{\rho-1}g(t){(^{\rho}I_{a+}^{\alpha}h)}(t)dt=\int_{a}^{b}t^{\rho-1}h(t){(^{\rho}I_{b-}^{\alpha}g)}(t)dt.
\end{gather}
\end{lemma}
\begin{proof}
The proof is straightforward. Using Dirichlet formula, we obtain
\begin{align*}
\int_{a}^{b}t^{\rho-1}g(t){(^{\rho}I_{a+}^{\alpha}h)(t)}dt&=\int_{a}^{b}t^{\rho-1}g(t)\int_{a}^{t}s^{\rho-1}{\bigg(\frac{t^\rho-s^\rho}{\rho}\bigg)}^{\alpha-1}\frac{h(s)}{\Gamma(\alpha)}{ds}~{dt}\\
&=\int_{a}^{b}s^{\rho-1}h(s)\int_{s}^{b}t^{\rho-1}{\bigg(\frac{t^\rho-s^\rho}{\rho}\bigg)}^{\alpha-1}\frac{g(t)}{\Gamma(\alpha)}{dt}~{ds}\\
&=\int_{a}^{b}t^{\rho-1}h(t){(^{\rho}I_{b-}^{\alpha}g)}(t)dt.
\end{align*}
\end{proof}
\begin{definition}\label{d8}\cite{oo}
The generalized Katugampola fractional derivatives, $^{\rho}D_{a+}^{\alpha,\beta}$ (left-sided) and $^{\rho}D_{b-}^{\alpha,\beta}$ (right-sided), of order $0<\alpha<1$ and type $0\leq\beta\leq1$ are respectively defined by
\begin{gather}\label{gkl}
{(^{\rho}D_{a+}^{\alpha,\beta}h)}(t)={({^{\rho}I_{a+}^{\beta(1-\alpha)}}\delta_\rho{{^{\rho}I_{a+}^{(1-\beta)(1-\alpha)}}}h)}(t),
\end{gather}
\begin{gather}\label{gkr}
{(^{\rho}D_{b-}^{\alpha,\beta}h)}(t)={(-{^{\rho}I_{b-}^{\beta(1-\alpha)}}\delta_\rho{{^{\rho}I_{b-}^{(1-\beta)(1-\alpha)}}}h)}(t)
\end{gather}
for the function to which right-hand side expressions exist and $\rho>0.$
\end{definition}
\begin{remark}\label{r1}
For $0<\alpha<1,0\leq\beta\leq1,$ the generalized Katugampola fractional derivative ${^{\rho}D_{a+}^{\alpha,\beta}}$ can be written in terms of Katugampola fractional derivative as
\begin{gather*}
  {^{\rho}D_{a+}^{\alpha,\beta}}={^{\rho}I_{a+}^{\beta(1-\alpha)}}\delta_\rho{^{\rho}I_{a+}^{1-\gamma}}={^{\rho}I_{a+}^{\beta(1-\alpha)}}~~{^{\rho}D_{a+}^{\gamma}},\qquad \gamma=\alpha+\beta(1-\alpha).
\end{gather*}
\end{remark}
\begin{lemma}\cite{ki}\label{l6}
Let $\alpha>0,0<a<b<\infty,0\leq\gamma<1$ and $g\in{C_{\gamma,\rho}[a,b]}.$ If $\alpha>\gamma,$ then
\begin{gather*}
  {(^{\rho}I_{a+}^{\alpha}g)}(a)=\lim_{x\to{a+}}{(^{\rho}I_{a+}^{\alpha}g)}(t)=0, \\
  {(^{\rho}I_{b-}^{\alpha}g)}(b)=\lim_{x\to{b-}}{(^{\rho}I_{b-}^{\alpha}g)}(t)=0.
\end{gather*}
\end{lemma}
\begin{theorem}\label{yi}[Young's inequality]
If $\varphi$ and $\eta$ are non-negative real numbers and $m$ and $m'$ are positive real numbers such that $\frac{1}{m}+\frac{1}{m'}=1,$ then we have
\begin{gather*}
\varphi\eta\leq\frac{\varphi^m}{m}+\frac{\eta^{m'}}{m'}.
\end{gather*}
\end{theorem}
\section{Non-existence result}
The proof of following theorem is based on the test function method developed by Mitidieri and Pokhazhaev in \cite{mp} and recently used in \cite{neh,nehh,newnec,newne,oldnex}.
\begin{theorem}
Assume that $\mu\in\mathbb{R}$ and $1<m<\frac{1+\mu}{1-\alpha},~\mu>-\alpha.$ Then, Problem \eqref{pex} does not admit global non-trivial solutions in $C_{1-\gamma,\rho}^{\gamma}[a,b],$ when $x_a>0$.
\end{theorem}
\begin{proof}
On the contrary, assume that a non-trivial solution exists for all time $t>a.$ Let $\phi\in{C^1([a,\infty))}$ be a test function satisfying $\phi\geq0$ and non-increasing such that
\begin{gather}\label{1}
\phi(t)=\begin{cases}
          1, &  a\leq t\leq\theta{T}, \\
          0, & t\geq{T},
        \end{cases}
\end{gather}
for some $T>a$ and some $\theta\leq\frac{1}{2}$ such that $a<\theta{T}<T$. Multiplying the inequality in \eqref{pex} by $\phi(t)$ and integrating over $[a,T],$ we obtain
\begin{gather}\label{2}
  \int_{a}^{T}\phi(t){(^{\rho}D_{a+}^{\alpha,\beta}x)}(t)dt\geq\int_{a}^{T}{\bigg(\frac{t^\rho-a^\rho}{\rho}\bigg)}^{\mu}{|x(t)|}^m\phi(t)dt.
\end{gather}
Observe that the integral in left-hand side exists and the one in the right-hand side exists for $m<\frac{1+\mu}{1-\alpha}$ when $x\in{C_{1-\gamma,\rho}^{\gamma}[a,b]}.$ Moreover, from Definition \ref{d8}, we can write \eqref{2} as
\begin{gather*}
  \int_{a}^{T}\phi(t){({^{\rho}I_{a+}^{\beta(1-\alpha)}}\delta_\rho{^{\rho}I_{a+}^{1-\gamma}}x)}(t)dt\geq\int_{a}^{T}{\bigg(\frac{t^\rho-a^\rho}{\rho}\bigg)}^{\mu}{|x(t)|}^m\phi(t)dt
\end{gather*}
\begin{gather}\label{3}
   \int_{a}^{T}t^{\rho-1}\phi(t){({^{\rho}I_{a+}^{\beta(1-\alpha)}}\delta_\rho{^{\rho}I_{a+}^{1-\gamma}}x)}(t)t^{1-\rho}dt\geq\int_{a}^{T}{\bigg(\frac{t^\rho-a^\rho}{\rho}\bigg)}^{\mu}{|x(t)|}^m\phi(t)dt.
\end{gather}
In accordance with Lemma \ref{l5} (after extending by zero outside $[a,T]$), from \eqref{3} we may deduce that
\begin{gather}\label{4}
\int_{a}^{T}t^{\rho-1}{(\delta_\rho{^{\rho}I_{a+}^{1-\gamma}}x)}(t){({^{\rho}I_{T-}^{\beta(1-\alpha)}}\phi)}(t)t^{1-\rho}dt\geq\int_{a}^{T}{\bigg(\frac{t^\rho-a^\rho}{\rho}\bigg)}^{\mu}{|x(t)|}^m\phi(t)dt.
\end{gather}
An integration by parts in \eqref{4} yields
\begin{align}\label{5}
[{(^{\rho}I_{a+}^{1-\gamma}x)}(t){({^{\rho}I_{T-}^{\beta(1-\alpha)}}\phi)}(t)]{\big|}_{t=a}^{T}-\int_{a}^{T}&{({^{\rho}I_{a+}^{1-\gamma}}x)}(t)\frac{d}{dt}{({^{\rho}I_{T-}^{\beta(1-\alpha)}}\phi)}(t)dt\nonumber\\
&\geq\int_{a}^{T}{\bigg(\frac{t^\rho-a^\rho}{\rho}\bigg)}^{\mu}{|x(t)|}^m\phi(t)dt.
\end{align}
Using ${({^{\rho}I_{T-}^{\beta(1-\alpha)}}\phi)}(T)=0$ (Lemma \ref{l6}) and initial condition ${({^{\rho}I_{a+}^{1-\gamma}}x)}(a+)=x_a,$ we obtain
\begin{align}\label{6}
-x_a{({^{\rho}I_{T-}^{\beta(1-\alpha)}}\phi)}(a)&-\int_{a}^{T}{({^{\rho}I_{a+}^{1-\gamma}}x)}(t)\frac{d}{dt}{({^{\rho}I_{T-}^{\beta(1-\alpha)}}\phi)}(t)dt\nonumber\\
&~~~\geq\int_{a}^{T}{\bigg(\frac{t^\rho-a^\rho}{\rho}\bigg)}^{\mu}{|x(t)|}^m\phi(t)dt.
\end{align}
Multiplying by $\frac{t^{\rho-1}}{t^{\rho-1}}$ inside the integral in the left-hand side of expression \eqref{6}, we see that
\begin{gather}\label{7}
{\hspace{-0.2cm}}L=\int_{a}^{T}{({^{\rho}I_{a+}^{1-\gamma}}x)}(t)(-\delta_\rho){({^{\rho}I_{T-}^{\beta(1-\alpha)}}\phi)}(t)\frac{dt}{t^{\rho-1}}\geq\int_{a}^{T}{\bigg(\frac{t^\rho-a^\rho}{\rho}\bigg)}^{\mu}{|x(t)|}^m\phi(t)dt.
\end{gather}
From Definition \ref{d7} (for $n=1$), we have
\begin{gather}\label{8}
L=\int_{a}^{T}{({^{\rho}I_{a+}^{1-\gamma}}x)}(t){({^{\rho}D_{T-}^{1-\beta(1-\alpha)}}\phi)}(t)dt\geq\int_{a}^{T}{\bigg(\frac{t^\rho-a^\rho}{\rho}\bigg)}^{\mu}{|x(t)|}^m\phi(t)dt.
\end{gather}
and from Lemma \ref{l4}, we see that
\begin{align*}
  L=\int_{a}^{T}{({^{\rho}I_{a+}^{1-\gamma}}x)}(t)&\bigg[\frac{\phi(T)}{\Gamma(\beta(1-\alpha))}{\bigg(\frac{T^\rho-a^\rho}{\rho}\bigg)}^{\beta(1-\alpha)}\\
  &-\int_{t}^{T}s^{\rho-1}{\bigg(\frac{s^\rho-t^\rho}{\rho}\bigg)}^{\beta(1-\alpha)-1}\frac{\phi^{'}(s)}{\Gamma(\beta(1-\alpha))}ds\bigg]dt.
\end{align*}
Since $\phi(T)=0$ we get
\begin{gather}\label{9}
 L=-\int_{a}^{T}{({^{\rho}I_{a+}^{1-\gamma}}x)}(t){({^{\rho}I_{T-}^{\beta(1-\alpha)}}\delta_\rho\phi)}(t)dt\geq\int_{a}^{T}{\bigg(\frac{t^\rho-a^\rho}{\rho}\bigg)}^{\mu}{|x(t)|}^m\phi(t)dt.
\end{gather}
Therefore
\begin{gather*}
L=-\int_{a}^{T}t^{\rho-1}{({^{\rho}I_{a+}^{1-\gamma}}x)}(t){({^{\rho}I_{T-}^{\beta(1-\alpha)}}\delta_\rho\phi)}(t)\frac{dt}{t^{\rho-1}}\geq\int_{a}^{T}{\bigg(\frac{t^\rho-a^\rho}{\rho}\bigg)}^{\mu}{|x(t)|}^m\phi(t)dt.
\end{gather*}
Lemma \ref{l5} again allows us to write
\begin{gather*}
L=-\int_{a}^{T}\delta_\rho\phi(t){({^{\rho}I_{a+}^{\beta(1-\alpha)}}{({^{\rho}I_{a+}^{1-\gamma}}x)})}(t)dt\geq\int_{a}^{T}{\bigg(\frac{t^\rho-a^\rho}{\rho}\bigg)}^{\mu}{|x(t)|}^m\phi(t)dt,
\end{gather*}
and by Lemma \ref{l1} we obtain
\begin{gather}\label{10}
L=-\int_{a}^{T}\delta_\rho\phi(t){{({^{\rho}I_{a+}^{1-\alpha}}x)}}(t)dt\geq\int_{a}^{T}{\bigg(\frac{t^\rho-a^\rho}{\rho}\bigg)}^{\mu}{|x(t)|}^m\phi(t)dt.
\end{gather}
On the other hand
\begin{align}\label{11}
\int_{a}^{T}\delta_\rho\phi(t){{({^{\rho}I_{a+}^{1-\alpha}}x)}}(t)dt&=\int_{a}^{T}\delta_\rho\phi(t)\int_{a}^{t}{\bigg(\frac{t^\rho-s^\rho}{\rho}\bigg)}^{-\alpha}\frac{s^{\rho-1}x(s)}{\Gamma(1-\alpha)}ds~dt\nonumber\\
&\leq\int_{a}^{T}|\delta_\rho\phi(t)|\int_{a}^{t}{\bigg(\frac{t^\rho-s^\rho}{\rho}\bigg)}^{-\alpha}\frac{s^{\rho-1}|x(s)|}{\Gamma(1-\alpha)}ds~dt.
\end{align}
As $\phi$ is non-increasing, we have $\phi(s)\geq\phi(t)$ for all $t\geq{s}$ and $\frac{1}{\phi^{{1}/{m}}(s)}\leq\frac{1}{\phi^{{1}/{m}}(t)}$ for $m>1.$ Also it is clear that $\phi^{'}(t)=0,\,\,t\in[a,\theta{T}].$ Therefore
\begin{align}\label{12}
L&\leq{\int_{a}^{T}|\delta_\rho\phi(t)|\int_{a}^{t}\frac{s^{\rho-1}}{\Gamma(1-\alpha)}{\bigg(\frac{t^\rho-s^\rho}{\rho}\bigg)}^{-\alpha}\frac{|x(s)|{\phi^{{1}/{m}}(s)}}{{\phi^{{1}/{m}}(s)}}ds~dt}\nonumber\\
&\leq{\int_{\theta{T}}^{T}\frac{|\delta_\rho\phi(t)|}{\phi^{{1}/{m}}(t)}{\int_{a}^{t}\frac{s^{\rho-1}}{\Gamma(1-\alpha)}{\bigg(\frac{t^\rho-s^\rho}{\rho}\bigg)}^{-\alpha}{|x(s)|{\phi^{{1}/{m}}(s)}}}ds~dt}.
\end{align}
Definition \ref{d4} allows us to write
\begin{gather}\label{13}
L\leq\int_{\theta{T}}^{T}t^{\rho-1}\frac{|\delta_\rho\phi(t)|}{\phi^{{1}/{m}}(t)}{\big({^{\rho}I_{a+}^{1-\alpha}{|x|{\phi^{{1}/{m}}}}}\big)}(t)\frac{dt}{t^{\rho-1}}.
\end{gather}
Moreover, it is easy to see that $\frac{\delta_\rho\phi(t)}{\phi^{{1}/{m}}(t)}\in{L_p},$ for otherwise, we consider $\phi^{\lambda}(t)$ with some sufficiently large $\lambda.$ Thus, we can apply Lemma \ref{l5} to obtain
\begin{gather}\label{14}
L\leq\int_{\theta{T}}^{T}{|x(t)|{\phi^{{1}/{m}}}}(t){\bigg({^{\rho}I_{T-}^{1-\alpha}}\frac{|\delta_\rho\phi|}{\phi^{{1}/{m}}}\bigg)}(t)dt.
\end{gather}
Next, we multiply by ${\big(\frac{t^\rho-a^\rho}{\rho}\big)}^{\frac{\mu}{m}}{\big(\frac{t^\rho-a^\rho}{\rho}\big)}^{\frac{-\mu}{m}}$ inside the integral in the right-hand side of \eqref{14}:
\begin{gather}\label{15}
L\leq\int_{\theta{T}}^{T}{\bigg({^{\rho}I_{T-}^{1-\alpha}}\frac{|\delta_\rho\phi|}{\phi^{{1}/{m}}}\bigg)}(t){|x(t)|{\phi^{{1}/{m}}}}(t)\frac{{\big(\frac{t^\rho-a^\rho}{\rho}\big)}^{\frac{\mu}{m}}}{{\big(\frac{t^\rho-a^\rho}{\rho}\big)}^{\frac{\mu}{m}}}dt.
\end{gather}
For $-\alpha<\mu<0,$ we have ${\big(\frac{t^\rho-a^\rho}{\rho}\big)}^{-\frac{\mu}{m}}<{\big(\frac{T^\rho-a^\rho}{\rho}\big)}^{-\frac{\mu}{m}}$ because $t<T$ and for $\mu>0$ we obtain
${\big(\frac{t^\rho-a^\rho}{\rho}\big)}^{-\frac{\mu}{m}}<{\theta}^{\frac{\mu}{m}}{\big(\frac{T^\rho-a^\rho}{\rho}\big)}^{-\frac{\mu}{m}}$
because $t>\theta T.$ It follows that
\begin{gather*}
  {\bigg(\frac{t^\rho-a^\rho}{\rho}\bigg)}^{-\frac{\mu}{m}}<\max{\{1,{\theta}^{{\mu}/{m}}\}}{\bigg(\frac{T^\rho-a^\rho}{\rho}\bigg)}^{-\frac{\mu}{m}}.
\end{gather*}
Therefore,
\begin{gather}\label{16}
L\leq\max{\{1,{\theta}^{{\mu}/{m}}\}}{\bigg(\frac{T^\rho-a^\rho}{\rho}\bigg)}^{-\frac{\mu}{m}}\int_{\theta{T}}^{T}{\bigg({^{\rho}I_{T-}^{1-\alpha}}\frac{|\delta_\rho\phi|}{\phi^{{1}/{m}}}\bigg)}(t){{\bigg(\frac{t^\rho-a^\rho}{\rho}\bigg)}^{\frac{\mu}{m}}}{|x(t)|{\phi^{{1}/{m}}}}(t)dt.
\end{gather}
A simple application of Young's inequality (Theorem \ref{yi}) with $m$ and $m'$ such that $\frac{1}{m}+\frac{1}{m'}=1$ in the right-hand side of \eqref{16} yields
\begin{align}\label{17}
L&\leq\frac{1}{m}\int_{\theta{T}}^{T}{\bigg({|x(t)|{\phi^{{1}/{m}}}}(t){{\bigg(\frac{t^\rho-a^\rho}{\rho}\bigg)}^{\frac{\mu}{m}}}\bigg)}^{m}(t)dt\nonumber\\
&\hspace{0.3cm}+\frac{{\big(\frac{T^\rho-a^\rho}{\rho}\big)}^{\frac{-\mu{m'}}{m}}{(\max{\{1,{\theta}^{{\mu}/{m}}\}})}^{m'}}{m'}
\int_{\theta{T}}^{T}{\bigg({^{\rho}I_{T-}^{1-\alpha}}\frac{|\delta_\rho\phi|}{\phi^{{1}/{m}}}\bigg)}^{m'}(t)dt\nonumber\\
&\leq\frac{1}{m}\int_{a}^{T}{|x(t)|}^{m}{\phi}(t){{\bigg(\frac{t^\rho-a^\rho}{\rho}\bigg)}^{{\mu}}}dt\nonumber\\
&\hspace{0.3cm}+\frac{{\big(\frac{T^\rho-a^\rho}{\rho}\big)}^{\frac{-\mu{m'}}{m}}{(\max{\{1,{\theta}^{{\mu}/{m}}\}})}^{m'}}{m'}
\int_{\theta{T}}^{T}{\bigg({^{\rho}I_{T-}^{1-\alpha}}\frac{|\delta_\rho\phi|}{\phi^{{1}/{m}}}\bigg)}^{m'}(t)dt.
\end{align}
Clearly, from \eqref{7} and \eqref{17}, we see that
\begin{align}\label{18}
\frac{{\big(\frac{T^\rho-a^\rho}{\rho}\big)}^{\frac{-\mu{m'}}{m}}{(\max{\{1,{\theta}^{{\mu}/{m}}\}})}^{m'}}{m'}
&\int_{\theta{T}}^{T}{\bigg({^{\rho}I_{T-}^{1-\alpha}}\frac{|\delta_\rho\phi|}{\phi^{{1}/{m}}}\bigg)}^{m'}(t)dt\nonumber\\
\geq{\bigg(1-\frac{1}{m}\bigg)}&\int_{a}^{T}{|x(t)|}^{m}{\phi}(t){{\bigg(\frac{t^\rho-a^\rho}{\rho}\bigg)}^{{\mu}}}dt
\end{align}
or
\begin{align}\label{19}
{\bigg(\frac{1}{m'}\bigg)}\int_{a}^{T}{|x(t)|}^{m}{\phi}(t){{\bigg(\frac{t^\rho-a^\rho}{\rho}\bigg)}^{{\mu}}}dt&\leq\frac{{\big(\frac{T^\rho-a^\rho}{\rho}\big)}^{\frac{-\mu{m'}}{m}}{(\max{\{1,{\theta}^{\frac{\mu}{m}}\}})}^{m'}}{m'}\nonumber\\
&~~~~\times\int_{\theta{T}}^{T}{\bigg({^{\rho}I_{T-}^{1-\alpha}}\frac{|\delta_\rho\phi|}{\phi^{{1}/{m}}}\bigg)}^{m'}(t)dt.
\end{align}
Therefore by Definition \ref{d5}, we have
\begin{align}\label{20}
\int_{a}^{T}{\bigg(\frac{t^\rho-a^\rho}{\rho}\bigg)}^{\mu}{\phi(t)}&{|x(t)|}^{m}dt
\leq\frac{{\big(\frac{{(\theta{T})}^{\rho}-a^\rho}{\rho}\big)}^{\frac{-\mu{m'}}{m}}}{{(\Gamma(1-\alpha))}^{m'}}{(\max{\{1,{\theta}^{{\mu}/{m}}\}})}^{m'}\nonumber\\
&\times\int_{{\theta}T}^{T}{\bigg(\int_{t}^{T}{\bigg(\frac{s^\rho-t^\rho}{\rho}\bigg)}^{-\alpha}s^{\rho-1}{\frac{|\delta_\rho\phi(s)|}{\phi^{1/m}(s)}}ds\bigg)}^{m'}(t)dt
\end{align}
The change of variable $t=\sigma{T}$ in the right-hand side integral yields
\begin{align}\label{21}
\int_{a}^{T}{\bigg(\frac{t^\rho-a^\rho}{\rho}\bigg)}^{\mu}&{\phi(t)}{|x(t)|}^{m}dt
\leq\frac{{\big(\frac{{(\theta{T})}^{\rho}-a^\rho}{\rho}\big)}^{\frac{-\mu{m'}}{m}}}{{(\Gamma(1-\alpha))}^{m'}}{(\max{\{1,{\theta}^{{\mu}/{m}}\}})}^{m'}\nonumber\\
\times\int_{{\theta}}^{1}&{\bigg(\int_{\sigma{T}}^{T}{\bigg(\frac{s^\rho-{(\sigma{T})}^\rho}{\rho}\bigg)}^{-\alpha}s^{\rho-1}{\frac{|\delta_\rho\phi(s)|}{\phi^{1/m}(s)}}ds\bigg)}^{m'}(\sigma)Td\sigma.
\end{align}
Another change of variable $s=rT$ therein yields
\begin{align}\label{22}
\int_{a}^{T}{\bigg(\frac{t^\rho-a^\rho}{\rho}\bigg)}^{\mu}&{\phi(t)}{|x(t)|}^{m}dt
\leq\frac{{\big(\frac{{(\theta{T})}^{\rho}-a^\rho}{\rho}\big)}^{\frac{-\mu{m'}}{m}}}{{(\Gamma(1-\alpha))}^{m'}}{(\max{\{1,{\theta}^{\frac{\mu}{m}}\}})}^{m'}\nonumber\\
\times\int_{{\theta}}^{1}&{\bigg(\int_{\sigma}^{T}{\bigg(\frac{{(rT)}^\rho-{(\sigma{T})}^\rho}{\rho}\bigg)}^{-\alpha}{(rT)}^{\rho-1}{\frac{|\delta_\rho\phi(r)|}{\phi^{1/m}(r)}}dr\bigg)}^{m'}(\sigma)Td\sigma.
\end{align}
We may assume that the integral in the right-hand side of \eqref{22} is bounded, i.e.
\begin{gather}\label{23}
\frac{{(\max{\{1,{\theta}^{{\mu}/{m}}\}})}^{m'}}{{(\Gamma(1-\alpha))}^{m'}}\int_{{\theta}}^{1}{\bigg(\int_{\sigma}^{T}{\bigg(\frac{{r}^\rho-{\sigma}^\rho}{\rho}\bigg)}^{-\alpha}{r}^{\rho-1}{\frac{|\delta_\rho\phi(r)|}{\phi^{1/m}(r)}}dr\bigg)}^{m'}(\sigma)d\sigma\leq{C},
\end{gather}
for some positive constant $C$, for otherwise we consider $\phi^{\lambda}(r)$ with some sufficiently large $\lambda$. Therefore
\begin{gather}\label{24}
\int_{a}^{T}{\bigg(\frac{t^\rho-a^\rho}{\rho}\bigg)}^{\mu}{\phi(t)}{|x(t)|}^{m}dt
\leq{C}{\bigg(\frac{{(\theta{T})}^{\rho}-a^\rho}{\rho}\bigg)}^{\frac{-\mu{m'}}{m}}.
\end{gather}
If $\mu>0,$ then
\begin{gather*}
{\bigg(\frac{{(\theta{T})}^{\rho}-a^\rho}{\rho}\bigg)}^{\frac{-\mu{m'}}{m}}\to0
\end{gather*}
as $T\to\infty.$ Finally, from \eqref{24}, we obtain
\begin{gather}\label{25}
  \lim_{T\to\infty}\int_{a}^{T}{\bigg(\frac{t^\rho-a^\rho}{\rho}\bigg)}^{\mu}{\phi(t)}{|x(t)|}^{m}dt=0.
\end{gather}
This is a contradiction since the solution is assumed to be nontrivial.\newline
In case $\mu=0$ we have $-\frac{\mu{m'}}{m}=0$ as the relation \eqref{24} ensures that
\begin{gather}\label{26}
  \lim_{T\to\infty}\int_{a}^{T}{\bigg(\frac{t^\rho-a^\rho}{\rho}\bigg)}^{\mu}{\phi(t)}{|x(t)|}^{m}dt\leq{C}.
\end{gather}
Moreover, it is clear that
\begin{align}\label{27}
{\bigg(\frac{{(\theta{T})}^\rho-a^\rho}{\rho}\bigg)}^{-\frac{\mu}{m}}
\int_{\theta{T}}^{T}&{\bigg({^{\rho}I_{T-}^{1-\alpha}}{\frac{|\delta_\rho\phi|}{\phi^{1/m}}}\bigg)}(t){\bigg(\frac{t^\rho-a^\rho}{\rho}\bigg)}^{\frac{\mu}{m}}{|x(t)|\phi^{1/m}(t)}dt\nonumber\\
&\leq{\bigg(\frac{{(\theta{T})}^\rho-a^\rho}{\rho}\bigg)}^{-\frac{\mu}{m}}{\bigg[\int_{\theta{T}}^{T}{\bigg({^{\rho}I_{T-}^{1-\alpha}}{\frac{|\delta_\rho\phi|}{\phi^{1/m}}}\bigg)}(t)dt\bigg]}^{\frac{1}{m'}}\nonumber\\
&\hspace{1cm}\times{\bigg[\int_{\theta{T}}^{T}{\bigg(\frac{t^\rho-a^\rho}{\rho}\bigg)}^{\frac{\mu}{m}}{|x(t)|\phi^{1/m}(t)}dt\bigg]}^{\frac{1}{m}}.
\end{align}
This relation \eqref{27} together with \eqref{16} implies that
\begin{gather}\label{28}
\int_{a}^{T}{\bigg(\frac{t^\rho-a^\rho}{\rho}\bigg)}^\mu\phi(t){|x(t)|}^mdt\leq{K}{\bigg[\int_{\theta{T}}^{T}{\bigg(\frac{t^\rho-a^\rho}{\rho}\bigg)}^{\frac{\mu}{m}}{|x(t)|\phi^{1/m}(t)}dt\bigg]}^{\frac{1}{m}}
\end{gather}
for some positive constant $K$, with
\begin{gather}\label{29}
\lim_{T\to\infty}\int_{\theta{T}}^{T}{\bigg(\frac{t^\rho-a^\rho}{\rho}\bigg)}^{\mu}{\phi(t)}{|x(t)|}^{m}dt=0
\end{gather}
due to the convergence of the integral in \eqref{26}. This is again contradiction.

If $\mu<0,$ we have ${\big(\frac{t^\rho-a^\rho}{\rho}\big)}^{-\frac{\mu}{m}}\leq{\big(\frac{T^\rho-a^\rho}{\rho}\big)}^{-\frac{\mu}{m}},$ because $-\frac{\mu}{m}>0$ and $t<T$. Similar to previous one, the expression $\frac{|\phi^{'}(r)|}{\phi^{1/m}(r)}$ may be assumed bounded and hence it can be shown that
\begin{gather}\label{30}
\int_{a}^{T}{\bigg(\frac{t^\rho-a^\rho}{\rho}\bigg)}^{\mu}{\phi(t)}{|x(t)|}^{m}dt\leq{C}{\bigg(\frac{T^\rho-a^\rho}{\rho}\bigg)}^{-m'-\mu{m'/m}}
\end{gather}
for some positive constant $C.$ This completes the proof.
\end{proof}
\section{Illustrations}
In this section, we prove the sharpness of exponent $\frac{\mu+1}{1-\alpha}$ i.e. the solutions exist for exponents strictly bigger than $\frac{\mu+1}{1-\alpha}.$ We need the following Lemma.
\begin{lemma}\label{l7}
Suppose that $0<\alpha<1,0\leq\beta\leq1$ and $\rho>0,\xi>0.$ Then,
the following identity holds for generalized Katugampola fractional derivative ${(^{\rho}D_{a+}^{\alpha,\beta})}:$
\begin{gather*}
\bigg({^{\rho}D_{a+}^{\alpha,\beta}}{\bigg(\frac{s^\rho-a^\rho}{\rho}\bigg)}^{\xi-1}\bigg)(t)=\frac{\Gamma(\xi)}{\Gamma(\xi-\alpha)}{\bigg(\frac{t^\rho-a^\rho}{\rho}\bigg)}^{\xi-\alpha-1},\quad t>a>0.
\end{gather*}
\end{lemma}
\begin{proof}
From Lemma \ref{l2}, for $\gamma=\alpha+\beta-\alpha\beta$ we see that
\begin{gather}\label{31}
\bigg({^{\rho}I_{a+}^{1-\gamma}}{\bigg(\frac{s^\rho-a^\rho}{\rho}\bigg)}^{\xi-1}\bigg)(t)=\frac{\Gamma(\xi)}{\Gamma(\xi-\gamma+1)}{\bigg(\frac{t^\rho-a^\rho}{\rho}\bigg)}^{\xi-\gamma}.
\end{gather}
Therefore
\begin{gather}\label{32}
\delta_\rho\bigg({^{\rho}I_{a+}^{1-\gamma}}{\bigg(\frac{s^\rho-a^\rho}{\rho}\bigg)}^{\xi-1}\bigg)(t)=
  \frac{\Gamma(\xi)}{\Gamma(\xi-\gamma)}{\bigg(\frac{t^\rho-a^\rho}{\rho}\bigg)}^{\xi-\gamma-1}
\end{gather}
which yields
\begin{gather}\label{33}
\bigg({^{\rho}D_{a+}^{\gamma}}{\bigg(\frac{s^\rho-a^\rho}{\rho}\bigg)}^{\xi-1}\bigg)(t)
  =\frac{\Gamma(\xi)}{\Gamma(\xi-\gamma)}{\bigg(\frac{t^\rho-a^\rho}{\rho}\bigg)}^{\xi-\gamma-1}
\end{gather}
Applying ${^{\rho}I_{a+}^{\beta(1-\alpha)}}$ on both sides of \eqref{33} and using Lemma \ref{l2}, again in the light of Definition \ref{d6} we conclude the Lemma.
\end{proof}
\begin{example}
Consider the following differential equation involving generalized Katugampola fractional derivative of order $0<\alpha<1$ and type $0\leq\beta\leq1:$
\begin{gather}\label{34}
{(^{\rho}D_{a+}^{\alpha,\beta}g)}(t)=\lambda{\bigg(\frac{t^\rho-a^\rho}{\rho}\bigg)}^{\mu}[g(t)]^{m},\quad t>a>0,m>1,
\end{gather}
with $\lambda,\mu\in\mathbb{R^+},\lambda\neq0,\rho>0.$
\end{example}
Look for the solution of the form $g(t)=c{\big(\frac{t^\rho-a^\rho}{\rho}\big)}^{\nu}$ for some $\nu\in\mathbb{R}.$\newline
Our aim is to find the values of $c$ and $\nu.$ By using Lemma \ref{l7}, we have
\begin{gather}\label{35}
{\bigg({^{\rho}D_{a+}^{\alpha,\beta}\bigg[c{\bigg(\frac{s^\rho-a^\rho}{\rho}\bigg)}^{\nu}\bigg]}\bigg)}(t)=\frac{c\Gamma(\nu+1)}{\Gamma(\nu-\alpha+1)}{\bigg(\frac{t^\rho-a^\rho}{\rho}\bigg)}^{\nu-\alpha}.
\end{gather}
Therefore, plugging this expression in \eqref{34}, we get
\begin{gather}\label{36}
\frac{c\Gamma(\nu+1)}{\Gamma(\nu-\alpha+1)}{\bigg(\frac{t^\rho-a^\rho}{\rho}\bigg)}^{\nu-\alpha}=\lambda{\bigg(\frac{t^\rho-a^\rho}{\rho}\bigg)}^{\mu}
{\bigg[c{\bigg(\frac{t^\rho-a^\rho}{\rho}\bigg)}^{\nu}\bigg]}^{m}.
\end{gather}
We obtain $\nu=\frac{\alpha+\mu}{1-m}$ and $c=\frac{\Gamma(\frac{\alpha+\mu}{1-m}+1)}{\lambda{\Gamma(\frac{m\alpha+\mu}{1-m}+1)}}.$ \newline
If $\frac{m\alpha+\mu}{1-m}>-1$ means $m>(\frac{1+\mu}{1-\alpha}),$ then \eqref{34} has the exact solution:
\begin{gather}\label{37}
  y(t)={\bigg[\frac{\Gamma(\frac{\alpha+\mu}{1-m}+1)}{\lambda{\Gamma(\frac{m\alpha+\mu}{1-m}+1)}}\bigg]}^{\frac{1}{m-1}}{\bigg(\frac{t^\rho-a^\rho}{\rho}\bigg)}^{\frac{\alpha+\mu}{1-m}}.
\end{gather}
This solution satisfies the initial condition with
\begin{gather*}
x_a={\bigg[\frac{\Gamma(\frac{\alpha+\mu}{1-m}+1)}{\lambda{\Gamma(\frac{m\alpha+\mu}{1-m}+1)}}\bigg]}^{\frac{1}{m-1}},\,\,\text{when}\,\,\frac{\alpha+\mu}{1-m}\geq\gamma-1>-1.
\end{gather*}
\begin{remark}
From the above non-existence theorem, we observe that no solutions can exists for certain values of $m$ and the exponent $\mu$ for all time $t>a.$ Clearly, we determine the range of values of $m$ for which solutions do not exists globally. Also note that, sufficient conditions for non-existence of solutions provide the necessary conditions for existence of solutions.
\end{remark}
\section{Conclusions} We successfully employed the fractional integration by parts and the method of test functions to obtain the necessary conditions for non-existence of global solutions for a wide class of fractional differential equations. The obtained results are well illustrated with suitable examples. Result discussed in this paper necessarily generalized / improved the existing results in literature.

\end{document}